\documentclass[<opyions>]{elsarticle}
\usepackage{amsfonts}
\usepackage{tipa}
\usepackage{mathrsfs}
\usepackage{amssymb}
\usepackage{latexsym,bm}
\usepackage{amsthm}
\usepackage{yhmath}
\usepackage{booktabs}
\usepackage{multirow}
\usepackage[center]{caption}
\usepackage{graphics}
\usepackage{subfig}
\usepackage{color}
\usepackage{algorithmic}
\usepackage{algorithm}
\usepackage{geometry}

\journal{arXiv.org}

\newtheorem{example}{Example}[section]
\newtheorem{remark}{Remark}[section]

\begin{document}
\begin{frontmatter}
\title{An unstructured mesh control volume method for two-dimensional space fractional diffusion equations with variable coefficients on convex domains}%\tnoteref{t1}} \tnotetext[t1]{This
%research has been supported by the } %\tnotetext[t2]{}
\author[els]{Libo~Feng}
%\ead{fenglibo2012@126.com}
%\author[rvt,rvh]{P.~Zhuang}
%\ead{zxy1104@xmu.edu.cn}
\author[els]{Fawang Liu\corref{cor1}}
\ead{f.liu@qut.edu.au}\cortext[cor1]{Corresponding author.}
\author[els,qv]{Ian Turner}
%\ead{i.turner@qut.edu.au}
\address[els]{School of Mathematical Sciences, Queensland University of
Technology,  GPO Box 2434, Brisbane, Qld 4001, Australia}
\address[qv]{Australian Research Council Centre of Excellence for Mathematical and Statistical Frontiers (ACEMS), Queensland University of Technology (QUT), Brisbane, Australia.}

\begin{abstract}
In this paper, we propose a novel unstructured mesh control volume method to deal with the space fractional derivative on arbitrarily shaped convex domains, which to the best of our knowledge is a new contribution to the literature. Firstly, we present the finite volume scheme for the two-dimensional space fractional diffusion equation with variable coefficients and provide the full implementation details for the case where the background interpolation mesh is based on triangular elements. Secondly, we explore the property of the stiffness matrix generated by the integral of space fractional derivative. We find that the stiffness matrix is sparse and not regular. Therefore, we choose a suitable sparse storage format for the stiffness matrix and develop a fast iterative method to solve the linear system, which is more efficient than using the Gaussian elimination method. Finally, we present several examples to verify our method, in which we make a comparison of our method with the finite element method for solving a Riesz space fractional diffusion equation on a circular domain. The numerical results demonstrate that our method can reduce CPU time significantly while retaining the same accuracy and approximation property as the finite element method. The numerical results also illustrate that our method is effective and reliable and can be applied to problems on arbitrarily shaped convex domains.
\end{abstract}

\end{frontmatter}
\section{Introduction}
In the past two decades, fractional differential equations have been applied in many fields of science \cite{a1,a2,a3,a4,a5,a6,a0}, in which space fractional diffusion equations are used to model the anomalous transport of solute in groundwater hydrology \cite{M1,Liu04}. For space fractional diffusion equations with constant coefficients, analytical solutions can be obtained by utilising the Fourier transform methods. However, many practical problems involve variable coefficients \cite{Bar,Che}, in which the diffusion velocity can vary over the solution domain. The work involving space fractional diffusion equations with variable coefficients is numerous. Meerschaert et al. \cite{M1,M3} considered the finite difference method for the one-dimensional one-sided and two-sided space fractional diffusion equations with variable coefficients, respectively. Zhang et al. \cite{Zhang} explored the homogeneous space-fractional advection-dispersion equation with space-dependent coefficients. Ding et al. \cite{Xiao} presented the weighted finite difference methods for a class of space fractional partial differential equations with variable coefficients.  Moroney and Yang \cite{Tim1,Tim2} proposed some fast preconditioners for the numerical solution of a class of two-sided nonlinear space-fractional diffusion equations with variable coefficients. Chen and Deng \cite{Den1} discussed the alternating direction implicit method to solve a two-dimensional two-sided space fractional convection-diffusion equation on a finite domain. Wang and Zhang \cite{Wang1} developed a high-accuracy preserving spectral Galerkin method for the Dirichlet boundary-value problem of a one-sided variable-coefficient conservative fractional diffusion equation. Feng et al. \cite{Feng1} proposed the finite volume method for a two-sided space-fractional diffusion equation with variable coefficients. Chen et al. \cite{Chen2} considered an inverse problem for identifying the fractional derivative indices in a two-dimensional space-fractional nonlocal model with variable diffusivity coefficients. Jia and Wang \cite{Wang2} presented a fast finite volume method for conservative space-fractional diffusion equations with variable coefficients. In \cite{Feng2}, Feng et al. presented a new second order finite difference scheme for a two-sided space-fractional diffusion equation with variable coefficients.

In fact, many mathematical models and problems from science and engineering must be computed on irregular domains and therefore seeking effective numerical methods to solve these problems on such domains is important. Although existing numerical methods for fractional diffusion equations are numerous \cite{a7,a8,a9,a10,a11,a12,a13,a14,a15,a16,a17,a18}, most of them are limited to regular domains and uniform meshes. Research involving unstructured meshes and irregular domains is sparse.   Yang et al. \cite{Yang14} proposed the finite volume scheme for a two-dimensional space-fractional reaction-diffusion equation based on the fractional Laplacian operator  $-(-\nabla^2)^{\frac{\alpha}{2}}$, which was computed using unstructured triangular meshes on a unit disk. Burrage et al. \cite{kevin} developed some
techniques for solving fractional-in-space reaction diffusion equations using the finite element method on both structured and unstructured grids.  Qiu et al. \cite{Qiu15} developed the nodal discontinuous Galerkin method for fractional diffusion equations on a two-dimensional domain with triangular meshes. Liu et al. \cite{Liu15} presented the semi-alternating direction method for a two-dimensional fractional FitzHugh-Nagumo monodomain model on an approximate irregular domain. Qin et al. \cite{Qin17} also used the implicit alternating direction method to solve a two-dimensional fractional Bloch-Torrey equation using an approximate irregular domain. Karaa et al. \cite{Karaa} proposed a finite volume element method implemented on an unstructured mesh for approximating the anomalous subdiffusion equations with a temporal fractional derivative. Yang et al. \cite{Yang17} established the unstructured mesh finite element method for the nonlinear Riesz space fractional diffusion equations on irregular convex domains. Fan et al. \cite{Fan} extended the unstructured mesh finite element method developed by Yang et al. \cite{Yang17} to the time-space fractional wave equation. Feng et al. \cite{Feng3} investigated the unstructured mesh finite element method for a two-dimensional time-space Riesz fractional diffusion equation on irregular arbitrarily shaped convex domains and a multiply-connected domain. Le et al. \cite{Le} studied the finite element approximation for a time-fractional diffusion problem on a domain with a re-entrant corner. To the best of our knowledge, the control volume finite element method (see Carr et al. \cite{Ian} for an illustration of the method applied to wood drying) has not been generalised to allow the solution of space fractional diffusion equations with variable coefficients.

In this paper, we will consider the unstructured mesh control volume method for the following two-dimensional space fractional diffusion equation with variable coefficients (2D SFDE-VC) \cite{Chen2} on an arbitrarily shaped convex domain:
\begin{align}
\frac{\partial u(x,y,t)}{\partial t}=&\frac{\partial }{\partial x}\bigg[K_1(x,y,t)\frac{\partial^{\alpha} u(x,y,t)}{\partial
x^{{\alpha}}}-K_2(x,y,t)\frac{\partial^{\alpha} u(x,y,t)}{\partial (-x)^{{\alpha}}}  \bigg]\nonumber\\
+&\frac{\partial}{\partial y}\bigg[K_3(x,y,t)\frac{\partial^{\beta} u(x,y,t)}{\partial
y^{{\beta}}}-K_4(x,y,t)\frac{\partial^{\beta} u(x,y,t)}{\partial (-y)^{{\beta}}}\bigg]\nonumber\\\label{e01}
+&f(x,y,t),\quad(x,y,t)\in \Omega\times (0,T],
\end{align}
subject to the initial condition
\begin{equation}\label{e02}
u(x,y,0)=\phi (x,y), \quad (x,y)\in\overline{\Omega},
\end{equation}
and boundary conditions
\begin{equation}\label{e03}
u(x,y,t)=0, \quad (x,y,t)\in \partial\Omega\times [0,T],
\end{equation}
where $0<\alpha,~\beta<1$, $K_i(x,y,t)\geq0$, $i=1,2,3,4$, $f(x,y,t)$ and $\phi(x,y)$ are assumed to be two known smooth functions. When the solution domain is rectangular $\Omega =(a,b)\times(c,d)$, we define the Riemman-Liouville fractional derivative as \cite{Pod99}:
\begin{align*}
\frac{\partial^{\alpha} u(x,y,t)}{\partial
x^\alpha}&={_{a}{D}_x^{\alpha}u(x,y,t)}=\frac{1}{\Gamma (1-\alpha )}\frac{\partial}{\partial x}\int_{a}^{x}{(x-s })^{-\alpha}u(s,y,t)~{d}s ,\\
\frac{\partial^{\alpha} u(x,y,t)}{\partial (-x)^{\alpha}}&={_{x}{D}_{b}^{\alpha}}u(x,y,t)=\frac{-1}{\Gamma (1-\alpha)}
\frac{\partial} {\partial x}\int_{x}^{b}{(s -x})^{-\alpha}u(s,y,t)~{d}s,\\
\frac{\partial^{\beta} u(x,y,t)}{\partial
y^\beta}&={_{c}{D}_y^{\beta}u(x,y,t)}=\frac{1}{\Gamma (1-\beta )}\frac{\partial}{\partial y}\int_{c}^{y}{(y-s })^{-\beta}u(x,s,t)~{d}s ,\\
\frac{\partial^{\beta} u(x,y,t)}{\partial (-y)^{\beta}}&={_{y}{D}_{d}^{\beta}}u(x,y,t)=\frac{-1}{\Gamma (1-\beta)}
\frac{\partial} {\partial y}\int_{y}^{d}{(s-y})^{-\beta}u(x,s,t)~{d}s.
\end{align*}

\begin{figure}[h]
\begin{center}
\scalebox{0.4}[0.4]{\includegraphics{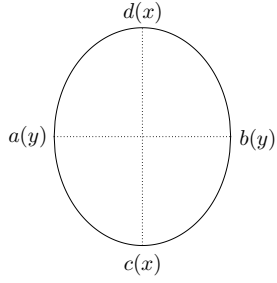}}
\caption{The illustration of a solution domain with curved boundary}
\label{fig1}
\end{center}
\end{figure}

When the boundary of the solution domain is nonconstant or curved, for example a convex domain shown in Figure \ref{fig1} with left boundary $a(y)$, right boundary $b(y)$, lower boundary $c(x)$ and upper boundary $d(x)$, we define the Riemman-Liouville fractional derivative as \cite{Feng3}:
\begin{align*}
\frac{\partial^{\alpha} u(x,y,t)}{\partial
x^\alpha}&={_{a(y)}{D}_x^{\alpha}u(x,y,t)}=\frac{1}{\Gamma (1-\alpha )}\frac{\partial}{\partial x}\int_{a(y)}^{x}{(x-s })^{-\alpha}u(s,y,t)~{d}s ,\\
\frac{\partial^{\alpha} u(x,y,t)}{\partial (-x)^{\alpha}}&={_{x}{D}_{b(y)}^{\alpha}}u(x,y,t)=\frac{-1}{\Gamma (1-\alpha)}
\frac{\partial} {\partial x}\int_{x}^{b(y)}{(s -x})^{-\alpha}u(s,y,t)~{d}s,\\
\frac{\partial^{\beta} u(x,y,t)}{\partial
y^\beta}&={_{c(x)}{D}_y^{\beta}u(x,y,t)}=\frac{1}{\Gamma (1-\beta )}\frac{\partial}{\partial y}\int_{c(x)}^{y}{(y-s })^{-\beta}u(x,s,t)~{d}s ,\\
\frac{\partial^{\beta} u(x,y,t)}{\partial (-y)^{\beta}}&={_{y}{D}_{d(x)}^{\beta}}u(x,y,t)=\frac{-1}{\Gamma (1-\beta)}
\frac{\partial} {\partial y}\int_{y}^{d(x)}{(s-y})^{-\beta}u(x,s,t)~{d}s.
\end{align*}
\begin{remark}\label{rem1.1}
When $K_i(x,y,t)~i=1,2,3,4$ take the special form
\begin{align*}
K_1(x,y,t)&=K_2(x,y,t)=-\frac{K_x}{2\cos\frac{\pi(1+\alpha)}{2}},\\
K_3(x,y,t)&=K_4(x,y,t)=-\frac{K_y}{2\cos\frac{\pi(1+\beta)}{2}},
\end{align*}
equation (\ref{e01}) can be written as the following Riesz space fractional diffusion equation \cite{Liu15,Yang17}
\begin{equation}\label{e04}
\frac{\partial u(x,y,t)}{\partial t}=K_x\frac{\partial^{1+\alpha}u(x,y,t)}{\partial
|x|^{1+\alpha}}+K_y\frac{\partial^{1+\beta}u(x,y,t)}{\partial
|y|^{1+\beta}}+f(x,y,t),
\end{equation}
where
\begin{align*}
\frac{\partial^{1+\alpha}u(x,y,t)}{\partial |x|^{1+\alpha}}&=-\frac{1}{2\cos\frac{\pi(1+\alpha)}{2}}\bigg[\frac{\partial^{1+\alpha} u(x,y,t)}{\partial
x^\alpha}+ \frac{\partial^{1+\alpha} u(x,y,t)}{\partial (-x)^{\alpha}}\bigg],\\
\frac{\partial^{1+\beta}u(x,y,t)}{\partial |y|^{1+\beta}}&=-\frac{1}{2\cos\frac{\pi(1+\beta)}{2}}\bigg[\frac{\partial^{1+\beta} u(x,y,t)}{\partial
y^\beta}+\frac{\partial^{1+\beta} u(x,y,t)}{\partial (-y)^{\beta}} \bigg].
\end{align*}
One important application of equation (\ref{e04}) is in the study of cardiac arrhythmias. In two dimensions, the fractional FitzHugh-Nagumo monodomain model can be rewritten as a two-dimensional Riesz space fractional reaction-diffusion model, which can be used to describe the propagation of the electrical potential in heterogeneous cardiac tissue \cite{Liu15,Bueno}. This electrophysiological model of the heart can describe how electrical currents flow through the heart controlling its contraction and can be used to ascertain the effects of certain drugs designed to treat heart problems.
\end{remark}

The major contribution of this paper is as follows.
\begin{itemize}
\item [$\bullet$] Different from \cite{Yang14} and \cite{Karaa}, we consider the control volume method for the two-dimensional space fractional diffusion equation with variable coefficients, in which the space fractional operator is either the Riemman-Liouville fractional derivative or Riesz space fractional derivative. To the best of our knowledge, this is a new contribution to the literature.
\item [$\bullet$] We propose a novel technique utilizing the control volume method implemented with an unstructured triangular mesh to deal with the space fractional derivative on an irregular convex domain, which we believe provides a very flexible solution strategy because our considered solution domain can be arbitrarily convex. Compared to the finite difference method in \cite{Liu15,Qin17},
our method requires fewer grid nodes to generate the meshes in the solution domain partition.
\item [$\bullet$] For the methods considered in this paper, we construct the control volumes using triangular meshes and transform the problem (\ref{e01}) from the solution domain to a single control volume. Then we integrate problem (\ref{e01}) over an arbitrary control volume and change the control volume integral to a line integral over the control volume faces, which is approximated by the midpoint approximation. Moreover, we utilise the linear basis function to approximate the fractional derivatives at the midpoints of the control volume faces, in which some numerical techniques are used to handle the non-locality of the fractional derivative of the basis function.
\item [$\bullet$] We explore the property of the stiffness matrix generated by the integral of space fractional derivative. We find that the stiffness matrix is sparse and not regular. Specially, the more small the maximum edge of the triangulation is, the more sparse of the stiffness matrix becomes. Therefore, we choose a suitable sparse storage format for the stiffness matrix and utilise the bi-conjugate gradient stabilized method (Bi-CGSTAB) iterative method to solve the linear system, which is more efficient than using the Gaussian elimination method.
\item [$\bullet$] We present several examples to verify our method, in which we make a comparison of our method with the finite element method proposed in \cite{Yang17} for solving the Riesz space fractional diffusion equation (\ref{e04}) on a circular domain. In \cite{Yang17}, the authors develop an algorithm to form the stiffness matrix and derive the bilinear operator as
    \begin{align*}
     A(u,v)&=\frac{K_x}{2\cos\frac{\pi(1+\alpha)}{2}}\Big\{ \Big({_{a(y)}{D}_x^{\frac{(1+\alpha)}{2}}u},{_{x}{D}_{b(y)}^{\frac{(1+\alpha)}{2}}v}\Big)
     +\Big({_{x}{D}_{b(y)}^{\frac{(1+\alpha)}{2}}u},{_{a(y)}{D}_x^{\frac{(1+\alpha)}{2}}v}\Big)\Big\}\\
     &+\frac{{K_y}}{2\cos\frac{\pi(1+\beta)}{2}}\Big\{ \Big({_{c(x)}{D}_y^{\frac{(1+\beta)}{2}}u},{_{y}{D}_{d(x)}^{\frac{(1+\beta)}{2}}v}\Big)
     +\Big({_{y}{D}_{d(x)}^{\frac{(1+\beta)}{2}}u},{_{c(x)}{D}_y^{\frac{(1+\beta)}{2}}v}\Big)\Big\}.
     \end{align*}
     The bilinear form involves 8 fractional derivative terms and the approximation of two-fold multiple integrals, which are approximated by Gauss quadrature. While for the control volume method, we use the following form to generate the stiffness matrix form,
     \begin{align*}
      &\frac{K_x}{2\cos\frac{\pi(1+\alpha)}{2}}\oint_{\Gamma_i}\bigg[\frac{\partial^{\alpha} u(x,y,t)}{\partial
      x^{{\alpha}}}-\frac{\partial^{\alpha} u(x,y,t)}{\partial (-x)^{{\alpha}}}  \bigg]{d}{y}\\
      -&\frac{{K_y}}{2\cos\frac{\pi(1+\beta)}{2}}\oint_{\Gamma_i}\bigg[\frac{\partial^{\beta} u(x,y,t)}{\partial
      y^{{\beta}}}-\frac{\partial^{\beta} u(x,y,t)}{\partial (-y)^{{\beta}}}\bigg]{d}{x},
     \end{align*}
     in which we only need to calculate 4 fractional derivative terms and the approximation of line integrals. The numerical results demonstrate that our method can reduce CPU time significantly while retaining the same accuracy and approximation property as the finite element method. The numerical results also illustrate that our method is effective and reliable and can be applied to problems on arbitrarily convex domains.
\end{itemize}

The outline of this paper is as follows. In Section \ref{sec2}, the unstructured mesh control volume method for problem (\ref{e01}) is proposed and the full implementation details are provided. Then the property of the stiffness matrix is explored and a fast iterative solver is developed for the linear system. In Section \ref{sec3}, several numerical examples are presented to verify the effectiveness of the method and comparisons are made with existing methods to highlight its computational performance. Finally, some conclusions of the work are drawn.

\section{Control volume finite element method}\label{sec2}

\begin{figure}[h]
\begin{center}
\scalebox{0.25}[0.25]{\includegraphics{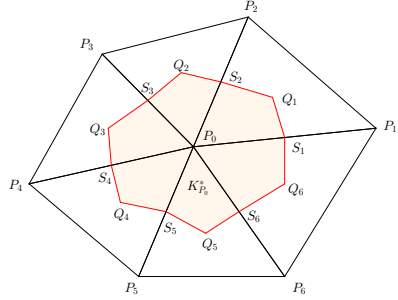}}
\caption{The illustration of a control volume}
\label{fig2}
\end{center}
\end{figure}

In this section, we will generalise the control volume method to solve equation (\ref{e01}), placing particular emphasis on the way
the Riemman-Liouville fractional derivatives are discretised in space. Firstly, we divide the solution domain $\Omega$ into a number of regular triangular regions. Let $\mathcal{T}_h$ denote this triangulation and $h$ be the maximum diameter of the triangular elements. Then we introduce the control volumes, which are constructed as follows. Let $M_h$ be a set of vertice,
\begin{eqnarray*}
M_h=\{P_i:~P_i~\text{is a vertex of the element}~ K\in \mathcal{T}_h ~\text{and}~ P_i\in \Omega\},
\end{eqnarray*}
and $M_h^0$ be the set of interior nodes in $\mathcal{T}_h$. We denote $P_0$ as the interior node of the triangulation $\mathcal{T}_h$ and $P_i~(i=1,2,\cdots,m)$ as its adjacent nodes (see Figure \ref{fig2} with $m=6$). Let $S_i~(i=1,2,\cdots,m)$ be the midpoints of the line segments $\overline{P_0P_i}$ and $Q_i~(i=1,2,\cdots,m)$ the barycenters of the triangle $\Delta P_0P_iP_{i+1}$ with $P_{m+1}=P_1$. The control volume $K_{P_0}^*$ is constructed by joining successively $S_1,~Q_1,~\cdots,~S_m,~Q_m,~S_1$ (see Figure \ref{fig2}). We call the line segments $\overline{S_iQ_i}$ and $\overline{Q_iS_{i+1}}$~($i=1,2,\cdots,m$ and $S_{m+1}=S_1)$ control volume faces. Consequently, each of the triangular elements is divided into three sub-domains by these control surfaces. These quadrilateral shapes are called sub-control volumes and are illustrated in Figure \ref{fig2} (for example, the quadrilateral $S_1Q_1S_2P_0$). Thus, a control volume consists of the sum of all neighbouring sub-control volumes that surround the given node $P_0$. The control volume is polygonal in shape and can be assembled in a straightforward and efficient manner at the element level. The flow across each control surface must be determined by an integral. Therefore, the finite volume method discretization process is initiated by utilising the integrated form of equation (\ref{e01}).

\begin{figure}[h]
\begin{center}
\scalebox{0.25}[0.25]{\includegraphics{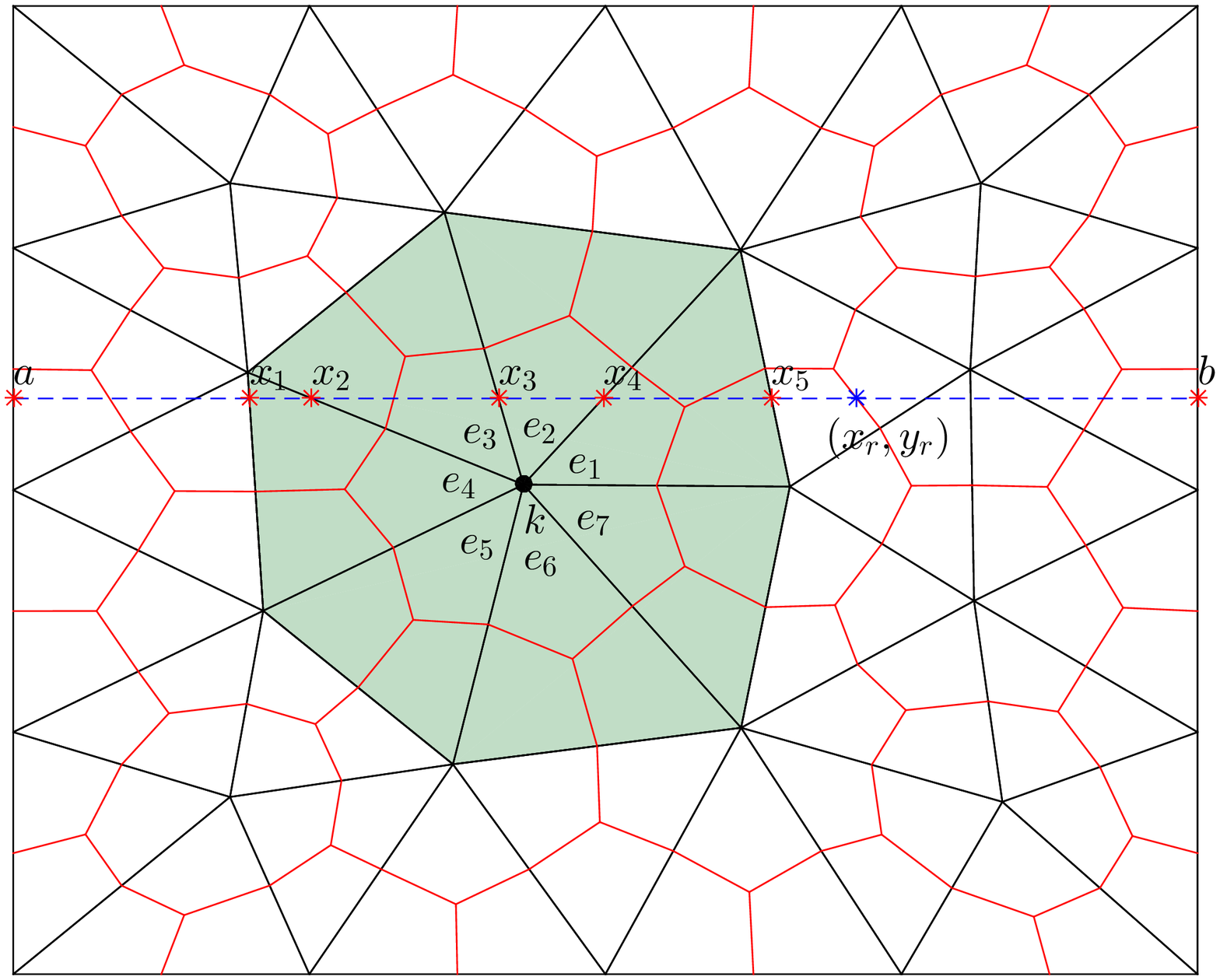}}
\caption{The illustration of line $y=y_r$ intersecting $n_q$ points with the support domain $\Omega_{e_k}$ of $l_k(x,y)$, where $(x_r,y_r)$ locates out of $\Omega_{e_k}$}\label{fig5}
\end{center}
\end{figure}

Integrating (\ref{e01}) over an arbitrary control volume $V_i$ ($i=1,2,\cdots,N_p$), yields
\begin{align}
\int_{V_i}\frac{\partial u(x,y,t)}{\partial t}~{d}V_i=&\int_{V_i}\frac{\partial }{\partial x}\bigg[K_1(x,y,t)\frac{\partial^{\alpha} u(x,y,t)}{\partial
x^{{\alpha}}}-K_2(x,y,t)\frac{\partial^{\alpha} u(x,y,t)}{\partial (-x)^{{\alpha}}}  \bigg]~{d}V_i\nonumber\\
+&\int_{V_i}\frac{\partial}{\partial y}\bigg[K_3(x,y,t)\frac{\partial^{\beta} u(x,y,t)}{\partial
y^{{\beta}}}-K_4(x,y,t)\frac{\partial^{\beta} u(x,y,t)}{\partial (-y)^{{\beta}}}\bigg]~{d}V_i\nonumber\\ \label{e05}
+&\int_{V_i}f(x,y,t)~{d}V_i.
\end{align}
Utilising a lumped mass approach for the time derivative and source term and applying Green's theorem to the other two integrals terms, gives
\begin{align}
\Delta{V_i}\frac{\partial u(x,y,t)}{\partial t}\bigg|_{(x_i,y_i)}=&\oint_{\Gamma_i}\bigg[K_1(x,y,t)\frac{\partial^{\alpha} u(x,y,t)}{\partial
x^{{\alpha}}}-K_2(x,y,t)\frac{\partial^{\alpha} u(x,y,t)}{\partial (-x)^{{\alpha}}}  \bigg]{d}{y}\nonumber\\
-&\oint_{\Gamma_i}\bigg[K_3(x,y,t)\frac{\partial^{\beta} u(x,y,t)}{\partial
y^{{\beta}}}-K_4(x,y,t)\frac{\partial^{\beta} u(x,y,t)}{\partial (-y)^{{\beta}}}\bigg]{d}{x}\nonumber\\\label{e06}
+&\Delta{V_i}f(x_i,y_i,t),
\end{align}
where $\Gamma_i$ is the boundary of control volume $V_i$. We assume the finite volume integration is an anticlockwise traversal and the outward unit normal surface vector to the control surface with $\Delta x=x_b-x_a$ and $\Delta y=y_b-y_a$.
%\begin{figure}[h]
%\begin{center}
%\scalebox{0.28}[0.28]{\includegraphics{F3_Vec1.eps}}
%\caption{A control volume face and the outward normal unit vector}\label{fig3}
%\end{center}
%\end{figure}
Denote $\Delta V_i$ and $\Delta V_{ij}$ the area of the control volume and the sub-control volume surrounding the point $(x_i,y_i)$, then we have
\begin{equation*}
\Delta V_i=\sum_{j=1}^{m_i}\Delta V_{ij},
\end{equation*}
where $m_i$ is the total number of sub-control volumes that make up the control volume associated with the node $i$. The integral term on the right-hand side of equation (\ref{e01}) is a line integral, which can be approximated by the midpoint approximation for each control surface. Hence, the first integral term in equation (\ref{e06}) can be rewritten as
\begin{align}
&\oint_{\Gamma_i}\bigg[K_1(x,y,t)\frac{\partial^{\alpha} u(x,y,t)}{\partial
x^{{\alpha}}}-K_2(x,y,t)\frac{\partial^{\alpha} u(x,y,t)}{\partial (-x)^{{\alpha}}}  \bigg]d{y}\nonumber\\\label{e07}
=&\sum_{j=1}^{m_i}\sum_{r=1}^{2}\bigg[K_1(x,y,t)\frac{\partial^{\alpha} u(x,y,t)}{\partial
x^{{\alpha}}}-K_2(x,y,t)\frac{\partial^{\alpha} u(x,y,t)}{\partial (-x)^{{\alpha}}}  \bigg]\bigg|_{(x_r,y_r)}\Delta y_{j,r}^i,
\end{align}
where $(x_r,y_r)$ is the mid-point of the control face (CF).
%\begin{figure}[h]
%\begin{center}
%\scalebox{0.25}[0.25]{\includegraphics{F4_CF.eps}}
%\caption{The illustration of control faces with mid-points}\label{fig4}
%\end{center}
%\end{figure}
Similarly, for the second integral term in equation (\ref{e06}), we have
\begin{align}
&\oint_{\Gamma_i}\bigg[K_3(x,y,t)\frac{\partial^{\beta} u(x,y,t)}{\partial
y^{{\beta}}}-K_4(x,y,t)\frac{\partial^{\beta} u(x,y,t)}{\partial (-y)^{{\beta}}}\bigg]d{x}\nonumber\\\label{e08}
=&\sum_{j=1}^{m_i}\sum_{r=1}^{2}\bigg[K_3(x,y,t)\frac{\partial^{\beta} u(x,y,t)}{\partial
y^{{\beta}}}-K_4(x,y,t)\frac{\partial^{\beta} u(x,y,t)}{\partial (-y)^{{\beta}}}\bigg]\bigg|_{(x_r,y_r)}\Delta x_{j,r}^i.
\end{align}
Substituting equations (\ref{e07}) and (\ref{e08}) into (\ref{e06}), we obtain
\begin{align}
&\Delta{V_i}\frac{\partial u(x,y,t)}{\partial t}\bigg|_{(x_i,y_i)}\nonumber\\
=&\sum_{j=1}^{m_i}\sum_{r=1}^{2}\bigg[K_1(x,y,t)\frac{\partial^{\alpha} u(x,y,t)}{\partial x^{{\alpha}}}-K_2(x,y,t)\frac{\partial^{\alpha} u(x,y,t)}{\partial (-x)^{{\alpha}}}  \bigg]\bigg|_{(x_r,y_r)}\Delta y_{j,r}^i\nonumber\\
-&\sum_{j=1}^{m_i}\sum_{r=1}^{2}\bigg[K_3(x,y,t)\frac{\partial^{\beta} u(x,y,t)}{\partial
y^{{\beta}}}-K_4(x,y,t)\frac{\partial^{\beta} u(x,y,t)}{\partial (-y)^{{\beta}}}\bigg]\bigg|_{(x_r,y_r)}\Delta x_{j,r}^i\nonumber\\\label{e09}
+&\Delta{V_i}f(x_i,y_i,t).
\end{align}
To discretise the time derivative in equation (\ref{e09}) at $t=t_n$, we use the backward Euler difference scheme
\begin{eqnarray}\label{e10}
\frac{\partial u(x,y,t_n)}{\partial t}=\frac{u(x,y,t_n)-u(x,y,t_{n-1})}{\tau}+O(\tau).
\end{eqnarray}
In the following, we discuss the spatial discretisation of $u(x,y,t_n)$. We consider the computation process for piecewise linear polynomials on the triangular element $e_p$, $p=1,2,...,N_e$, where $N_e$ is the total number of triangles. Then, within element $e_p$, the field function $u^p(x,y)$ can be written as
\begin{eqnarray*}
u^p(x,y)=\sum_{j=1}^{3} u_{j}~\varphi_{j}(x,y)+O(h^2),
\end{eqnarray*}
where the triangle vertices are numbered in a counter-clockwise order as $1,~2,~3$ and the basis function $\varphi_{j}(x,y)$ is defined as
\begin{align*}
&\varphi_{j}(x,y)\Big|_{(x,y)\in e_p}=\frac{1}{2\Delta_{e_p}}(a_{j}\,x+b_{j}~y+c_{j}),\quad \varphi_{j}(x,y)\Big|_{(x,y)\notin e_p}=0,\\
&a_1=y_2-y_3,~a_2=y_3-y_1,~a_3=y_1-y_2,\\
&b_1=x_3-x_2,~b_2=x_1-x_3,~b_3=x_2-x_1,\\
&c_1=x_2y_3-x_3y_2,~c_2=x_3y_1-x_1y_3,~c_3=x_1y_2-x_2y_1,
\end{align*}
where $\Delta_{e_p}$ is the area of triangle element $p$. It is well-known that
\begin{eqnarray*}
\varphi_{j}(x_{i},y_{i})=\delta_{{i}{j}},\quad {i},~{j}=1,~2,~3,
\end{eqnarray*}
where $\delta$ is the Kronecker function. With these local field functions and basis functions, we can obtain a global approximation of $u(x,y)$ for the whole triangulation:
\begin{eqnarray*}
u(x,y)=\sum_{k=1}^{N_p} u_{k}~l_{k}(x,y)+O(h^2),
\end{eqnarray*}
where $l_{k}(x,y)$ is the new basis function whose support domain is $\Omega_{e_{k}}$ (see Figure \ref{fig5} the green polygonal domain) and $N_p$ is the total number of vertices on the convex domain $\Omega$.

Now, we denote $u_h(x,y,t_n)$ as the approximation solution of $u(x,y,t_n)$ and write $u_h(x,y,t_n)$ in the form
\begin{eqnarray}\label{e11}
u_h(x,y,t_n)=\sum_{k=1}^{N_p} u_{k}^n~l_{k}(x,y),
\end{eqnarray}
where $u_{k}^n$ are the coefficients that are to be solved for. Substituting equations (\ref{e10}) and (\ref{e11}) into equation (\ref{e09}), we discretise equation (\ref{e09}) at $t=t_n$ as follows:
\begin{align}
&\Delta V_i\sum_{k=1}^{N_p}\frac{u_{k}^n-u_{k}^{n-1}}{\tau}l_{k}(x_i,y_i)\nonumber\\
=&\sum_{k=1}^{N_p}\sum_{j=1}^{m_i}\sum_{r=1}^{2}u_k^n\bigg[K_1(x,y,t)\frac{\partial^{\alpha} l_k(x,y)}{\partial x^{{\alpha}}}-K_2(x,y,t)\frac{\partial^{\alpha} l_k(x,y)}{\partial (-x)^{{\alpha}}}  \bigg]\bigg|_{(x_r,y_r)}\Delta y_{j,r}^i\nonumber\\
-&\sum_{k=1}^{N_p}\sum_{j=1}^{m_i}\sum_{r=1}^{2}u_k^n\bigg[K_3(x,y,t)\frac{\partial^{\beta} l_k(x,y)}{\partial
y^{{\beta}}}-K_4(x,y,t)\frac{\partial^{\beta} l_k(x,y)}{\partial (-y)^{{\beta}}}\bigg]\bigg|_{(x_r,y_r)}\Delta x_{j,r}^i\nonumber\\\label{e12}
+&\Delta{V_i}f(x_i,y_i,t_n).
\end{align}
Using the fact that
\begin{equation*}
l_{k}(x_i,y_i)=\left\{\begin{array}{ll}
1,& i=k,\\
0,& i \neq k,
\end{array}\right.
\end{equation*}
we obtain
\begin{align}
&\Delta V_i\frac{u_{i}^n-u_{i}^{n-1}}{\tau}\nonumber\\
=&\sum_{k=1}^{N_p}\sum_{j=1}^{m_i}\sum_{r=1}^{2}u_k^n\bigg[K_1(x,y,t)\frac{\partial^{\alpha} l_k(x,y)}{\partial x^{{\alpha}}}-K_2(x,y,t)\frac{\partial^{\alpha} l_k(x,y)}{\partial (-x)^{{\alpha}}}  \bigg]\bigg|_{(x_r,y_r)}\Delta y_{j,r}^i\nonumber\\
-&\sum_{k=1}^{N_p}\sum_{j=1}^{m_i}\sum_{r=1}^{2}u_k^n\bigg[K_3(x,y,t)\frac{\partial^{\beta} l_k(x,y)}{\partial
y^{{\beta}}}-K_4(x,y,t)\frac{\partial^{\beta} l_k(x,y)}{\partial (-y)^{{\beta}}}\bigg]\bigg|_{(x_r,y_r)}\Delta x_{j,r}^i\nonumber\\\label{e13}
+&\Delta{V_i}f(x_i,y_i,t_n).
\end{align}
Equation (\ref{e13}) can be written in the following matrix form
\begin{eqnarray}\label{e14}
\mathbf{A}\frac{\mathbf{U}^n-\mathbf{U}^{n-1}}{\tau}=\mathbf{M}\mathbf{U}^n+\mathbf{A}\mathbf{F}^n,
\end{eqnarray}
where $\mathbf{A}=$diag$~[\Delta V_1, \Delta V_2, \ldots,\Delta V_{N_p}]$, $\mathbf{U}^n=[u_1^n,u_2^n,\ldots,u_{N_p}^n]^T$, $\mathbf{F}^n=[f(x_1,y_1,t_n),f(x_2,$ $y_2,t_n),\ldots,f(x_{N_p},y_{N_p},t_n)]^T$. Rearranging we obtain
\begin{eqnarray}\label{e15}
(\mathbf{A}-\tau \mathbf{M})\mathbf{U}^n=\mathbf{A}\mathbf{U}^{n-1}+\tau \mathbf{A}\mathbf{F}^n.
\end{eqnarray}
To form matrix $\mathbf{M}$, we need to calculate the fractional derivative of the basis function $l_k(x,y)$. In the following, we focus on the calculation of $\frac{\partial^{\alpha} l_k(x,y)}{\partial x^{{\alpha}}}$, $\frac{\partial^{\alpha} l_k(x,y)}{\partial (-x)^{{\alpha}}}$, $\frac{\partial^{\beta} l_k(x,y)}{\partial y^{{\beta}}}$ and $\frac{\partial^{\beta} l_k(x,y)}{\partial (-y)^{{\beta}}}$ at $(x_r,y_r)$. To evaluate $\frac{\partial^{\alpha} l_k(x,y)}{\partial x^{{\alpha}}}\big|_{(x_r,y_r)}$ and $\frac{\partial^{\alpha} l_k(x,y)}{\partial (-x)^{{\alpha}}}\big|_{(x_r,y_r)}$, suppose that line $y=y_r$ intersects $n_q$ points with the support domain $\Omega_{e_k}$ of $l_k(x,y)$ (see Figure \ref{fig5} with $n_q=5$).

\begin{figure}[h]
\begin{center}
\scalebox{0.25}[0.25]{\includegraphics{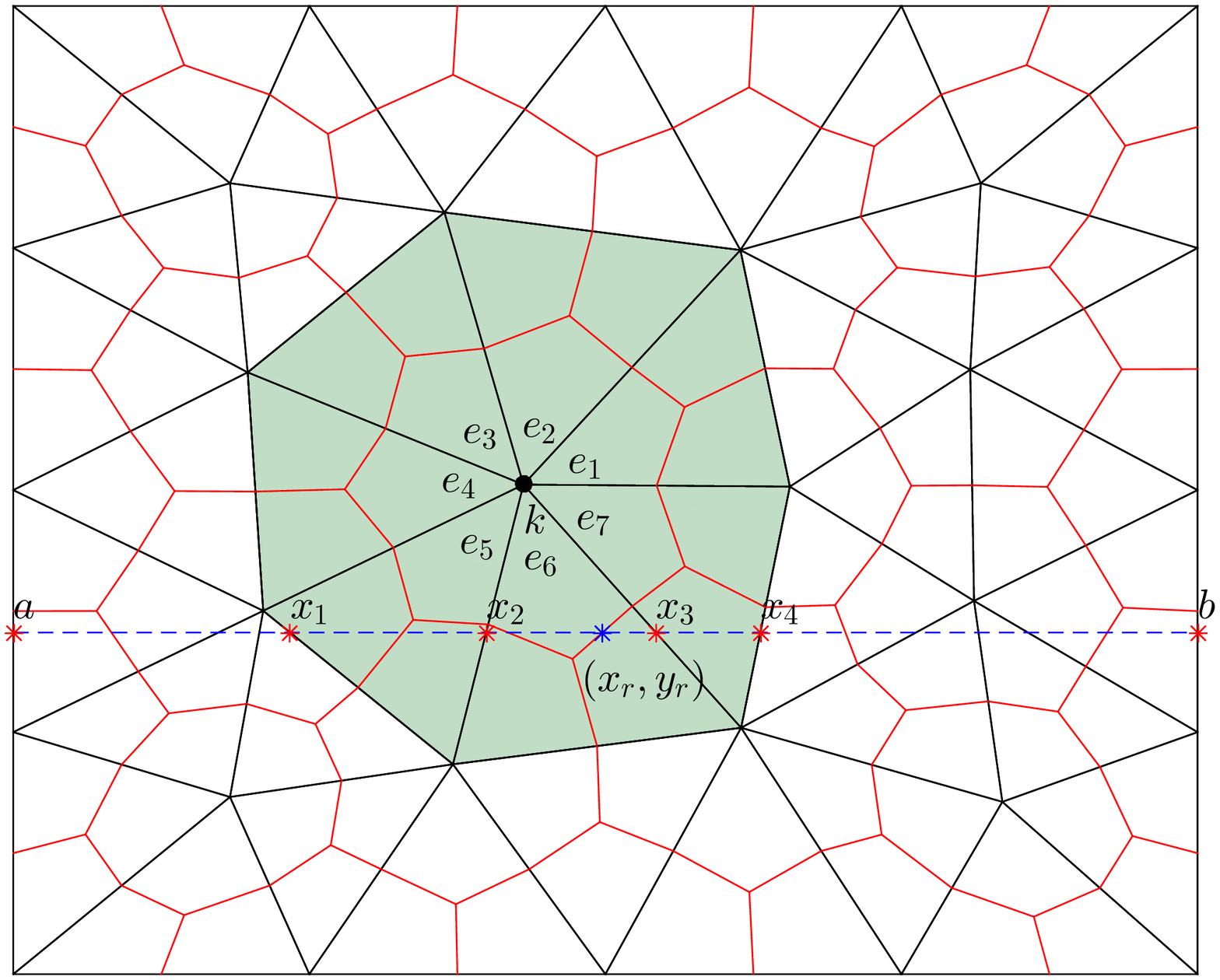}}
\caption{The illustration of line $y=y_r$ intersecting $n_q$ points with the support domain $\Omega_{e_k}$ of $l_k(x,y)$, where $(x_r,y_r)$ locates in $\Omega_{e_k}$}\label{fig6}
\end{center}
\end{figure}

Then we have
\begin{align*}
\frac{\partial^{\alpha} l_k(x,y)}{\partial x^{{\alpha}}}\bigg|_{(x_r,y_r)}&=\frac{\partial^{\alpha} l_k(x,y_r)}{\partial x^{{\alpha}}}\bigg|_{x=x_r},\\
\frac{\partial^{\alpha} l_k(x,y)}{\partial (-x)^{{\alpha}}}\bigg|_{(x_r,y_r)}&=\frac{\partial^{\alpha} l_k(x,y_r)}{\partial (-x)^{{\alpha}}}\bigg|_{x=x_r}.
\end{align*}
Using the important observation that
\begin{equation*}
l_k(x,y_r)=\left\{
\begin{array}{cc}
0, & a\leq x\leq x_1,\\
\varphi_{k4}(x,y_r),&x_1\leq x\leq x_2,\\
\varphi_{k3}(x,y_r),&x_2\leq x\leq x_3,\\
\varphi_{k2}(x,y_r),&x_3\leq x\leq x_4,\\
\varphi_{k1}(x,y_r),&x_4\leq x\leq x_5,\\
0, &x_5\leq x\leq b,\\
\end{array} \right.
\end{equation*}
where $\varphi_{kp}(x,y)$ is the basis function of node $k$ on the triangular element $e_p$, we obtain
\begin{align}
&\frac{\partial^{\alpha} l_k(x,y_r)}{\partial x^{{\alpha}}}\bigg|_{x=x_r}\nonumber\\
&=\bigg( \frac{1}{\Gamma(1-\alpha)}\frac{\partial}{\partial x} \int_{a}^{x} (x-\xi)^{-\alpha}l_k(\xi,{y_r})d\xi\bigg)\bigg|_{x=x_r}\nonumber\\
&=\bigg[ \frac{1}{\Gamma(1-\alpha)}\frac{\partial}{\partial x} \bigg(\int_{a}^{x_1}+\int_{x_1}^{x_2}+\int_{x_2}^{x_3}+\int_{x_3}^{x_4}+\int_{x_4}^{x_5}+\int_{x_5}^{x}\bigg) (x-\xi)^{-\alpha}l_k(\xi,{y_r})d\xi\bigg]\bigg|_{x=x_r}\nonumber\\\label{e16}
&=\bigg[ \frac{1}{\Gamma(1-\alpha)}\frac{\partial}{\partial x} \bigg(\int_{x_1}^{x_2}+\int_{x_2}^{x_3}+\int_{x_3}^{x_4}+\int_{x_4}^{x_5}\bigg) (x-\xi)^{-\alpha}l_k(\xi,{y_r})d\xi\bigg]\bigg|_{x=x_r}.
\end{align}
As $l_k(x,{y_r})$ is a linear function on each sub integral interval, equation (\ref{e16}) can be evaluated using integration by parts over each sub integral interval. For the right fractional derivative of $l_k(x,{y_r})$ at $(x_r,y_r)$, we obtain
\begin{align} \label{e17}
\frac{\partial^{\alpha} l_k(x,y_r)}{\partial (-x)^{{\alpha}}}\bigg|_{x=x_r}=\bigg( \frac{-1}{\Gamma(1-\alpha)}\frac{\partial}{\partial x} \int_{x}^{b} (\xi-x)^{-\alpha}l_k(\xi,{y_r})d\xi\bigg)\bigg|_{x=x_r}=0.
\end{align}
Now we consider the case that point $(x_r,y_r)$ is in the support domain $\Omega_{e_k}$ of $l_k(x,y)$. Suppose that line $y=y_r$ intersects $n_q$ points with the support domain $\Omega_{e_k}$  (see Figure \ref{fig6} with $n_q=4$). In this case, we have

\begin{equation*}
l_k(x,y_r)=\left\{
\begin{array}{cc}
0, & a\leq x\leq x_1,\\
\varphi_{k5}(x,y_r),&x_1\leq x\leq x_2,\\
\varphi_{k6}(x,y_r),&x_2\leq x\leq x_3,\\
\varphi_{k7}(x,y_r),&x_3\leq x\leq x_4,\\
0, &x_4\leq x\leq b.\\
\end{array} \right.
\end{equation*}
Then
\begin{align}
\frac{\partial^{\alpha} l_k(x,y_r)}{\partial x^{{\alpha}}}\bigg|_{x=x_r}&=\bigg( \frac{1}{\Gamma(1-\alpha)}\frac{\partial}{\partial x} \int_{a}^{x} (x-\xi)^{-\alpha}l_k(\xi,{y_r})d\xi\bigg)\bigg|_{x=x_r}\nonumber\\
&=\bigg[ \frac{1}{\Gamma(1-\alpha)}\frac{\partial}{\partial x} \bigg(\int_{a}^{x_1}+\int_{x_1}^{x_2}+\int_{x_2}^{x}\bigg) (x-\xi)^{-\alpha}l_k(\xi,{y_r})d\xi\bigg]\bigg|_{x=x_r}\nonumber\\\label{e18}
&=\bigg[ \frac{1}{\Gamma(1-\alpha)}\frac{\partial}{\partial x} \bigg(\int_{x_1}^{x_2}+\int_{x_2}^{x}\bigg) (x-\xi)^{-\alpha}l_k(\xi,{y_r})d\xi\bigg]\bigg|_{x=x_r},
\end{align}
and
\begin{align}
\frac{\partial^{\alpha} l_k(x,y_r)}{\partial (-x)^{{\alpha}}}\bigg|_{x=x_r}&=\bigg( \frac{-1}{\Gamma(1-\alpha)}\frac{\partial}{\partial x} \int_{x}^{b} (\xi-x)^{-\alpha}l_k(\xi,{y_r})d\xi\bigg)\bigg|_{x=x_r}\nonumber\\
&=\bigg[ \frac{-1}{\Gamma(1-\alpha)}\frac{\partial}{\partial x} \bigg(\int_{x}^{x_3}+\int_{x_3}^{x_4}+\int_{x_4}^{b}\bigg) (\xi-x)^{-\alpha}l_k(\xi,{y_r})d\xi\bigg]\bigg|_{x=x_r}\nonumber\\\label{e19}
&=\bigg[ \frac{-1}{\Gamma(1-\alpha)}\frac{\partial}{\partial x} \bigg(\int_{x}^{x_3}+\int_{x_3}^{x_4}\bigg) (\xi-x)^{-\alpha}l_k(\xi,{y_r})d\xi\bigg]\bigg|_{x=x_r}.
\end{align}
If line $y=y_r$ intersects zero points with the support domain $\Omega_{e_k}$, then we have
\begin{align} \label{e20}
\frac{\partial^{\alpha} l_k(x,y_r)}{\partial x^{{\alpha}}}\bigg|_{x=x_r}=0,\quad\frac{\partial^{\alpha} l_k(x,y_r)}{\partial (-x)^{{\alpha}}}\bigg|_{x=x_r}=0.
\end{align}
The calculation of $\frac{\partial^{\beta} l_k(x,y)}{\partial y^{{\beta}}}$ and $\frac{\partial^{\beta} l_k(x,y)}{\partial (-y)^{{\beta}}}$ at $(x_r,y_r)$ can be derived in a similar manner for the $y$ direction.  Finally, we summarise the whole computation process in the following algorithm (see Algorithm \ref{alg1}).

\begin{algorithm}[h]
\caption{Unstructured mesh CVM for solving 2D SFDE-VC }
\label{alg1}
\begin{algorithmic}[1]
\STATE{Partition the convex domain $\Omega$ with unstructured triangular elements $e_p$ and save the element information (node number, coordinates, and element number );}
\FOR{ $p=1,2,\cdots,N_e$}
\STATE{Find the barycenters of each triangular element $e_p$, form the control faces, sub-control volumes and save the sub-control volume information (the midpoint coordinates of each side of the triangular elements $e_p$, the midpoint coordinates $(x_r,y_r)$ of each control faces, etc.);}
\STATE{Calculate the areas of the sub-control volumes and control volumes, form matrix $\mathbf{A}$;}
\FOR{ $k=1,2,\cdots,N_p$}
\STATE{Find the support domain $\Omega_{e_k}$;}
\STATE{Find the points of intersection by $y={y_r}$ with $\Omega_{e_k}$ and calculate $\frac{\partial^{\alpha} l_k(x,y)}{\partial x^{{\alpha}}}\big|_{(x_r,y_r)}$,$\frac{\partial^{\alpha} l_k(x,y)}{\partial (-x)^{{\alpha}}}\big|_{(x_r,y_r)}$;}
\STATE{Find the points of intersection by $x={x_r}$ with $\Omega_{e_r}$ and calculate $\frac{\partial^{\beta} l_k(x,y)}{\partial y^{{\beta}}}\big|_{(x_r,y_r)}$,
$\frac{\partial^{\beta} l_k(x,y)}{\partial (-y)^{{\beta}}}\big|_{(x_r,y_r)}$;}
\ENDFOR
\STATE{Form the matrix $\mathbf{M}$;}
\STATE{Form the vector $\mathbf{F}^n$;}
\ENDFOR
\STATE{Solve the linear system (\ref{e15}) and obtain $\mathbf{U}^n$.}
\end{algorithmic}
\end{algorithm}

\begin{remark}\label{rem2.1}
When the boundary of the solution domain is nonconstant or curved, all of the above calculation is applicable as well.
\end{remark}
Here, we discuss the structure of matrix $\mathbf{M}$. Firstly, the matrix $\mathbf{M}$ generated by scheme (\ref{e13}) is sparse and not regular. Then we explore the sparsity of matrix $\mathbf{M}$ for different $h$. Table \ref{tab1} shows the size and density (nonzero entries percentage) of matrix $\mathbf{M}$ for different $h$ where we can observe that as $h$ decreases the density of matrix $\mathbf{M}$ reduces significantly. We can infer that when $h$ is small enough, matrix $\mathbf{M}$ is extremely sparse and this facilitates the use of a sparse matrix storage format to reduce the memory usage of our computational method. Furthermore, we employ an efficient sparse iterative solver Bi-CGSTAB \cite{Van} to solve the linear system (\ref{e15}) (see Algorithm \ref{alg2}), which is more efficient than using Gaussian elimination method. The CPU time comparison of the two methods is studied numerically in Example \ref{ex1}.
%\begin{figure}[h]
%\begin{center}
%\scalebox{0.25}[0.25]{\includegraphics{F7_h8.eps}}
%\caption{Sparsity pattern of matrix $\mathbf{M}$ for $h=1.6759\times 10^{-1}$. The size of $\mathbf{M}$ is  64$\times$64. Blue points indicate the nonzero entries}\label{fig7}
%\end{center}
%\end{figure}
\begin{table}[h]
\begin{center}
\caption{The size and density of matrix $\mathbf{M}$ for different $h$ on a square domain $[0,1]\times[0,1]$}\label{tab1}
\begin{tabular}{ccc}
\hline
\hline
       $h$   &  Size             & Density   \\
\hline
  5.2693E-01 &  4$\times$4         &   100\%       \\
  3.1123E-01 &  15$\times$15       &   86.667\%   \\
  1.6759E-01 &  64$\times$64       &   57.715\%   \\
  8.6682E-02 &  258$\times$258     &   34.002\%   \\
  4.3719E-02 &  1115$\times$1115   &   17.705\%	\\
  2.3063E-02 &  5255$\times$5255   &   8.517\%     \\
\hline
\hline
\end{tabular}
\end{center}
\end{table}

\begin{algorithm}[H]
\caption{The Bi-CGSTAB algorithm }
\label{alg2}
\begin{algorithmic}[1]
\STATE{Define $\mathbf{A}_0=\mathbf{A}-\tau \mathbf{M}$, use a sparse matrix storage format to store $\mathbf{A}_0$;}
\STATE{In each time level $t_n$, $\mathbf{x}_0=\mathbf{U}^{n-1}$, $\mathbf{b}=\mathbf{A}\mathbf{U}^{n-1}+\tau \mathbf{A} \mathbf{F}^{n}$;}
\STATE{Compute $\mathbf{r}_0=\mathbf{b}-\mathbf{A}_0\mathbf{x}_0$, $\hat{\mathbf{r}}_0$ is an arbitrary vector, such that $(\hat{\mathbf{r}}_0, \mathbf{r}_0)\neq0$. We choose $\hat{\mathbf{r}}_0=
\mathbf{r}_0$;}
\STATE{Let $\rho_0=\alpha_0=\omega_0=1$, $v_0=p_0=0$;}
\FOR{ $i=1,2,3,\cdots,$ }
\STATE{$\rho_i=(\hat{\mathbf{r}}_0, \mathbf{r}_{i-1});$}
\STATE{$\beta_0=(\rho_i/\rho_{i-1})(\alpha_{i-1}/\omega_{i-1})$;}
\STATE{$\mathbf{p}_i=\mathbf{r}_{i-1}+\beta_0(\mathbf{p}_{i-1}-\omega_{i-1}\mathbf{v}_{i-1})$;}
\STATE{$\mathbf{v}_{i}=\mathbf{A}_0\mathbf{p}_i$, $\alpha_{i}=\rho_i/(\hat{\mathbf{r}}_0, \mathbf{v}_i)$;}
\STATE{$\mathbf{s}=\mathbf{r}_{i-1}-\alpha_{i}\mathbf{v}_i$, $\mathbf{t}_0=\mathbf{A}_0\mathbf{s}$;}
\STATE{$\omega_{i}=(\mathbf{t}_0,\mathbf{s})/(\mathbf{t}_0,\mathbf{t}_0)$;}
\STATE{$\mathbf{x}_i=\mathbf{x}_{i-1}+\alpha_i \mathbf{p}_i+\omega_i \mathbf{s}$;}
\STATE{if $\mathbf{x}_i$ is accurate enough then quit;}
\STATE{$\mathbf{r}_i=\mathbf{s}-\omega_i \mathbf{t}_0$;}
\ENDFOR
\STATE{$\mathbf{U}^{n}=\mathbf{x}_i$.}
\end{algorithmic}
\end{algorithm}

\section{Discussion of Numerical Results}\label{sec3}

In this section, we provide some numerical examples to verify the effectiveness of our method presented in Section \ref{sec2}.  We adopt linear polynomials on triangles and define $h$ as the maximum length of the triangle edges.  $N_e$ is taken as the number of triangles in $\mathcal{T}_h$. Here, the numerical computations were carried out using MATLAB R2014b on a Dell desktop with configuration: Intel(R) Core(TM) i7-4790, 3.60 GHz and 16.0 GB RAM. We use the following formula to calculate the convergence order:
\begin{eqnarray*}
\mbox{Order}=\frac{\log({E(h_1)}/{E(h_2)})}{\log({h_1}/{h_2})}.
\end{eqnarray*}
\begin{example}\label{ex1}
Firstly, we consider the following 2D SFDE-VC on a rectangular domain
\begin{align*}
\frac{\partial u(x,y,t)}{\partial t}=&\frac{\partial }{\partial x}\bigg[K_1(x,y,t)\frac{\partial^{\alpha} u(x,y,t)}{\partial
x^{{\alpha}}}-K_2(x,y,t)\frac{\partial^{\alpha} u(x,y,t)}{\partial (-x)^{{\alpha}}}  \bigg]\nonumber\\
+&\frac{\partial}{\partial y}\bigg[K_3(x,y,t)\frac{\partial^{\beta} u(x,y,t)}{\partial
y^{{\beta}}}-K_4(x,y,t)\frac{\partial^{\beta} u(x,y,t)}{\partial (-y)^{{\beta}}}\bigg]\nonumber\\
+&f(x,y,t),\quad(x,y,t)\in \Omega\times (0,T],
\end{align*}
subject to
\begin{align*}
&u(x,y,0)=x^2(1-x)^2y^2(1-y)^2,\quad  (x,y)\in\overline{\Omega},\\
&u(x,y,t)=0,\quad  (x,y,t)\in \partial\Omega\times [0,T].
\end{align*}
where $\Omega=(0,1)\times(0,1)$, $T=1$,
\begin{align*}
f(x,y,t)&=2tx^2(1-x)^2y^2(1-y)^2-\Big[\frac{\partial K_1(x,y,t)}{\partial x}\cdot p(x,\alpha)+K_1(x,y,t)\cdot p(x,1+\alpha)\\
&-\frac{\partial K_2(x,y,t)}{\partial x}\cdot p(1-x,\alpha)+K_2(x,y,t)\cdot p(1-x,1+\alpha)\Big]y^2(1-y)^2(t^2+1)\\
&-\Big[\frac{\partial K_3(x,y,t)}{\partial y}\cdot p(y,\beta)+K_3(x,y,t)\cdot p(y,1+\beta)-\frac{\partial K_4(x,y,t)}{\partial y}\cdot p(1-y,\beta)\\
&+K_4(x,y,t)\cdot p(1-y,1+\beta)\Big]x^2(1-x)^2(t^2+1),\\
p(z,r)&=\frac{\Gamma(3)}{\Gamma(3-r)}z^{2-r}-\frac{2\Gamma(4)}{\Gamma(4-r)}z^{3-r}+\frac{\Gamma(5)}{\Gamma(5-r)}z^{4-r}.
\end{align*}
This is a two-dimensional anomalous diffusion model, which can describe anomalous transport in heterogeneous porous media and can be used
to explain the region-scale anomalous dispersion with heavy tails \cite{Chen2}.
\end{example}

The exact solution of this problem is given by $u(x,y,t)=(t^2+1)x^2(1-x)^2y^2(1-y)^2$.  Here, we consider three different coefficient cases \cite{Feng2}: linear coefficients $K_1(x,y,t)=2-x$, $K_2(x,y,t)=2+x$, $K_3(x,y,t)=2-y$, $K_4(x,y,t)=2+y$, quadratic coefficients $K_1(x,y,t)=2-x^2$, $K_2(x,y,t)=2+x^2$, $K_3(x,y,t)=2-y^2$, $K_4(x,y,t)=2+y^2$ and exponential coefficients $K_1(x,y,t)=3-e^x$, $K_2(x,y,t)=3+e^x$, $K_3(x,y,t)=3-e^y$, $K_4(x,y,t)=3+e^y$. The numerical results are given in Tables \ref{tab2} to \ref{tab4}. Table \ref{tab2} illustrates the $L_{2}$ error, $L_{\infty}$ error and corresponding convergence order of $h$ for the linear coefficient case for different $\alpha$, $\beta$ with $\tau=10^{-3}$ at $t=1$. Tables \ref{tab3} and \ref{tab4} show the $L_{2}$ error, $L_{\infty}$ error and corresponding convergence order of $h$ for the quadratic coefficient case and exponential coefficient case, respectively. From these tables we can see that the convergence order of both the $L_{2}$ error and $L_{\infty}$ error is $2-\max\{\alpha,\beta\}$ order \cite{Feng1} and the numerical results are in excellent agreement with the exact solution, which demonstrates the effectiveness of the numerical method. We can also observe that with $h$ deceasing, the CPU time grows considerably, which we believe is mainly due to the non-locality of the fractional derivative of the basis function and the computational cost to generate the matrix $\mathbf{M}$. In addition, we give a comparison between the Bi-CGSTAB and Gaussian elimination. In the Bi-CGSTAB solver, we set  $10^{-10}$ as the stopping criterion and the maximum iteration number is $10^{2}$. Table \ref{tab5} displays the consumed CPU time of these two algorithms at $t=1$ with $\tau=10^{-3}$, $\alpha=0.3$, $\beta=0.5$, $K_1(x,y,t)=2-x$, $K_2(x,y,t)=2+x$, $K_3(x,y,t)=2-y$, $K_4(x,y,t)=2+y$ for different $h$. Compared to the Gaussian elimination, Bi-CGSTAB has significantly reduced 90\% of the computational time for $h=4.3719E-02$. Another advantage of Bi-CGSTAB to be mentioned is that the average iteration number does not appear to increase significantly as $h$ decreases. Here, the average iteration number is approximately $10$ regardless of the model dimensions. We conclude that the Bi-CGSTAB solver is more efficient than  Gaussian elimination for solving this problem.
%\begin{figure}[h]
%\begin{center}
%\scalebox{0.35}[0.35]{\includegraphics{F8_h4.eps}}
%\scalebox{0.35}[0.35]{\includegraphics{F8_h8.eps}}\\
%\scalebox{0.35}[0.35]{\includegraphics{F8_h16.eps}}
%\scalebox{0.35}[0.35]{\includegraphics{F8_h32.eps}}
%\caption{The rectangular domain partitioned by unstructured meshes with control volumes for $h\approx 3.1123\times 10^{-1}, 1.6759\times 10^{-1}, 8.6682\times 10^{-2}, 4.3719\times 10^{-2}$, respectively}\label{fig8}
%\end{center}
%\end{figure}
\begin{table}[h]
\begin{center}
\caption{The $L_{2}$ error, $L_{\infty}$ error, convergence order and CPU time of $h$ with $\tau=10^{-3}$ for the linear coefficient case at $t=1$}\label{tab2}
\begin{tabular}{ccccccc}
\hline
\hline
             &     $h$ & $L_{2}$ error & Order & $L_{\infty}$ error & Order &Time\\
\hline
             & 3.1123E-01 & 3.5684E-04  & --   & 1.4774E-03 & --    & 4.90s\\

$\alpha=0.3$ & 1.6759E-01 & 1.0880E-04  & 1.92 & 4.3735E-04 & 1.97  & 19.50s\\

$\beta=0.5$  & 8.6682E-02 & 2.2391E-05  & 2.40 & 1.3895E-04 & 1.74  & 2.30min\\

             & 4.3719E-02 & 6.9379E-06  & 1.71 & 3.7632E-05 & 1.91  & 28.42min\\
\hline
             & 3.1123E-01 & 3.7935E-04  & --   & 1.4827E-03 & --    & 4.91s\\

$\alpha=0.4$ & 1.6759E-01 & 1.2435E-04  & 1.80 & 4.2971E-04 & 2.00  & 19.98s\\

$\beta=0.8$  & 8.6682E-02 & 2.5152E-05  & 2.42 & 1.3725E-04 & 1.73  & 2.36min\\

             & 4.3719E-02 & 7.2675E-06  & 1.81 & 3.5722E-05 & 1.97  & 28.56min\\
\hline
             & 3.1123E-01 & 3.9259E-04  & --   & 1.3844E-03 & --    & 4.91s\\

$\alpha=0.7$ & 1.6759E-01 & 1.4100E-04  & 1.65 & 4.1957E-04 & 1.93  & 19.87s\\

$\beta=0.9$  & 8.6682E-02 & 2.8670E-05  & 2.42 & 1.4117E-04 & 1.65  & 2.37min\\

             & 4.3719E-02 & 7.5385E-06  & 1.95 & 3.3666E-05 & 2.09  & 28.47min\\
\hline
\hline
\end{tabular}
\end{center}
\end{table}

\begin{table}[h]
\begin{center}
\caption{The $L_{2}$ error, $L_{\infty}$ error, convergence order and CPU time of $h$ with $\tau=10^{-3}$ for the quadratic coefficient case at $t=1$}\label{tab3}
\begin{tabular}{ccccccc}
\hline
\hline
             &     $h$ & $L_{2}$ error & Order & $L_{\infty}$ error & Order &Time\\
\hline
             & 3.1123E-01 & 3.1608E-04  & --   & 1.3430E-03 & --     & 4.97s\\

$\alpha=0.3$ & 1.6759E-01 & 1.0064E-04  & 1.85 & 4.0906E-04 & 1.92   & 20.48s\\

$\beta=0.5$  & 8.6682E-02 & 2.0661E-05  & 2.40 & 1.3852E-04 & 1.64   & 2.45min\\

             & 4.3719E-02 & 6.2709E-06  & 1.74 & 3.7584E-05 & 1.91   & 28.69min\\
\hline
             & 3.1123E-01 & 3.6299E-04  & --   & 1.4108E-03 & --     & 4.88s\\

$\alpha=0.4$ & 1.6759E-01 & 1.2145E-04  & 1.77 & 4.1614E-04 & 1.97   & 20.51s\\

$\beta=0.8$  & 8.6682E-02 & 2.4646E-05  & 2.42 & 1.3823E-04 & 1.67   & 2.46min\\

             & 4.3719E-02 & 6.7517E-06  & 1.89 & 3.3858E-05 & 2.06   & 28.78min\\
\hline
             & 3.1123E-01 & 3.8524E-04  & --   & 1.3424E-03 & --     & 4.97s\\

$\alpha=0.7$ & 1.6759E-01 & 1.3952E-04  & 1.64 & 4.0669E-04 & 1.93   & 20.56s\\

$\beta=0.9$  & 8.6682E-02 & 2.8522E-05  & 2.41 & 1.4126E-04 & 1.60   & 2.44min\\

             & 4.3719E-02 & 7.1520E-06  & 2.02 & 3.1880E-05 & 2.17   & 28.68min\\
\hline
\hline
\end{tabular}
\end{center}
\end{table}

\begin{table}[h]
\begin{center}
\caption{The $L_{2}$ error, $L_{\infty}$ error, convergence order and CPU time of $h$ with $\tau=10^{-3}$ for the exponential coefficient case at $t=1$}\label{tab4}
\begin{tabular}{ccccccc}
\hline
\hline
             &     $h$ & $L_{2}$ error & Order & $L_{\infty}$ error & Order &Time \\
\hline
             & 3.1123E-01 & 5.1809E-04  & --   & 1.9033E-03 & --    & 4.97s\\

$\alpha=0.3$ & 1.6759E-01 & 1.6296E-04  & 1.87 & 5.3973E-04 & 2.04  & 20.62s\\

$\beta=0.5$  & 8.6682E-02 & 3.8817E-05  & 2.18 & 1.6032E-04 & 1.84  & 2.45min\\

             & 4.3719E-02 & 1.1574E-05  & 1.77 & 4.8226E-05 & 1.76  & 28.46min\\
\hline
             & 3.1123E-01 & 4.5022E-04  & --   & 1.6750E-03 & --    & 4.93s\\

$\alpha=0.4$ & 1.6759E-01 & 1.4896E-04  & 1.79 & 1.0117E-04 & 2.01  & 20.52s\\

$\beta=0.8$  & 8.6682E-02 & 3.4126E-05  & 2.24 & 4.8309E-04 & 1.84  & 2.45min\\

             & 4.3719E-02 & 1.1238E-05  & 1.62 & 4.3016E-05 & 1.76  & 28.66min\\
\hline
             & 3.1123E-01 & 4.2412E-04  & --   & 1.4994E-03 & --    & 4.93s\\

$\alpha=0.7$ & 1.6759E-01 & 1.5286E-04  & 1.65 & 4.6520E-04 & 1.89  & 20.50s\\

$\beta=0.9$  & 8.6682E-02 & 3.3401E-05  & 2.31 & 1.4533E-04 & 1.76  & 2.45min\\

             & 4.3719E-02 & 1.0565E-05  & 1.68 & 4.0322E-05 & 1.87  & 28.56min\\
\hline
\hline
\end{tabular}
\end{center}
\end{table}

\begin{table}[H]
\begin{center}
\caption{Comparison of the consumed CPU time of Gaussian elimination versus Bi-CGSTAB}\label{tab5}
\begin{tabular}{cccc}
\hline
\hline
     $N_e$   &     $h$    & Gauss elimination & Bi-CGSTAB \\
\hline
     44      & 3.1123E-01 & 4.90s  & 4.90s\\

     158     & 1.6759E-01 & 22.57s  & 19.50s\\

     578     & 8.6682E-02 & 5.39min &2.30min\\

     2356    & 4.3719E-02 & 5.48h  & 28.42min\\
\hline
\hline
\end{tabular}
\end{center}
\end{table}

\begin{example}
Next, we consider the following two-dimensional Riesz space fractional diffusion equation on a circular domain, which can be used to describe the propagation of the electrical potential in heterogeneous cardiac tissue \cite{Liu15,Yang17,Bueno}.
\begin{equation}\label{e21}
\left\{\begin{array}{l}
\frac{\partial u(x,y,t)}{\partial t}=K_x\frac{{{\partial }^{{1+\alpha}}}u(x,y,t)}{\partial
{{|x|}^{{1+\alpha}}}}+K_y\frac{{{\partial }^{{1+\beta}}}u(x,y,t)}{\partial
{{|y|}^{1+\beta}}}+f(x,y,t),\quad (x,y,t)\in \Omega\times (0,T],\\
u(x,y,0)=({x^2}+{y^2}-1)^2, \quad (x,y)\in\overline{\Omega},\\
u(x,y,t)=0,\quad (x,y,t)\in \partial\Omega\times [0,T],
\end{array}\right.
\end{equation}
where $\Omega=\{(x,y)|{x^2}+{y^2}<1\}$, $K_x=1$, $K_y=1$, $T=1$,
\begin{align*}
f(x,y,t)&=-e^{-t}({x^2}+{y^2}-1)^2\\
&+\frac{e^{-t}}{2\cos((1+\alpha)/2\pi)}\bigg[\Big(f_1(x,a_0,\alpha)+g_1(x,b_0,\alpha)\Big)+({2y^2}
-{2})\Big(f_2(x,a_0,\alpha)+g_2(x,b_0,\alpha)\Big)\\
&+(y^2-1)^2\Big(f_3(x,a_0,\alpha)+g_3(x,b_0,\alpha)\Big)\bigg]\\
&+\frac{e^{-t}}{2\cos((1+\beta)/2\pi)}\bigg[\Big(f_1(y,c_0,\beta)+g_1(y,d_0,\beta)\Big)+({2x^2}
-{2})\Big(f_2(y,c_0,\beta)+g_2(y,d_0,\beta)\Big)\\
&+(x^2-1)^2\Big(f_3(y,c_0,\beta)+g_3(y,d_0,\beta)\Big)\bigg],\\
&a_0=-\sqrt{1-y^2}, ~ b_0=\sqrt{1-y^2},~c_0=-\sqrt{1-x^2}, ~ d_0=\sqrt{1-x^2},\\
&f_1(x,a,\alpha)={_aD_x^{1+\alpha}(x^4)},~f_2(x,a,\alpha)={_aD_x^{1+\alpha}(x^2)},~f_3(x,a,\alpha)={_aD_x^{1+\alpha}(1)},\\
&g_1(x,b,\alpha)={_xD_b^{1+\alpha}(x^4)},~g_2(x,b,\alpha)={_xD_b^{1+\alpha}(x^2)},~g_3(x,b,\alpha)={_xD_b^{1+\alpha}(1)}.
\end{align*}
\end{example}

\begin{figure}[h]
\begin{center}
%\scalebox{0.3}[0.3]{\includegraphics{F9_Circular_4.eps}}
%\scalebox{0.3}[0.3]{\includegraphics{F9_Circular_8.eps}}\\
\scalebox{0.3}[0.3]{\includegraphics{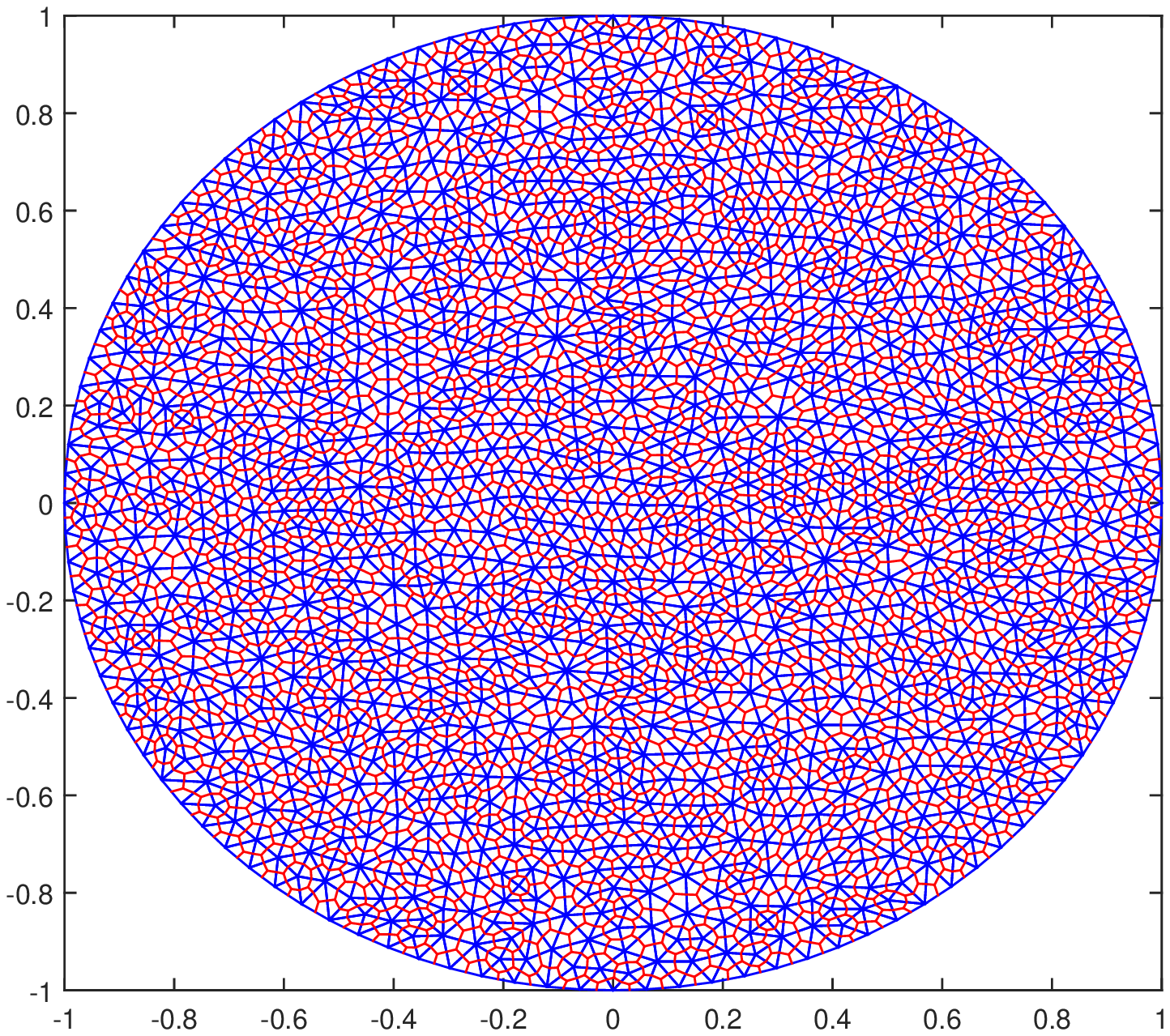}}
\scalebox{0.3}[0.3]{\includegraphics{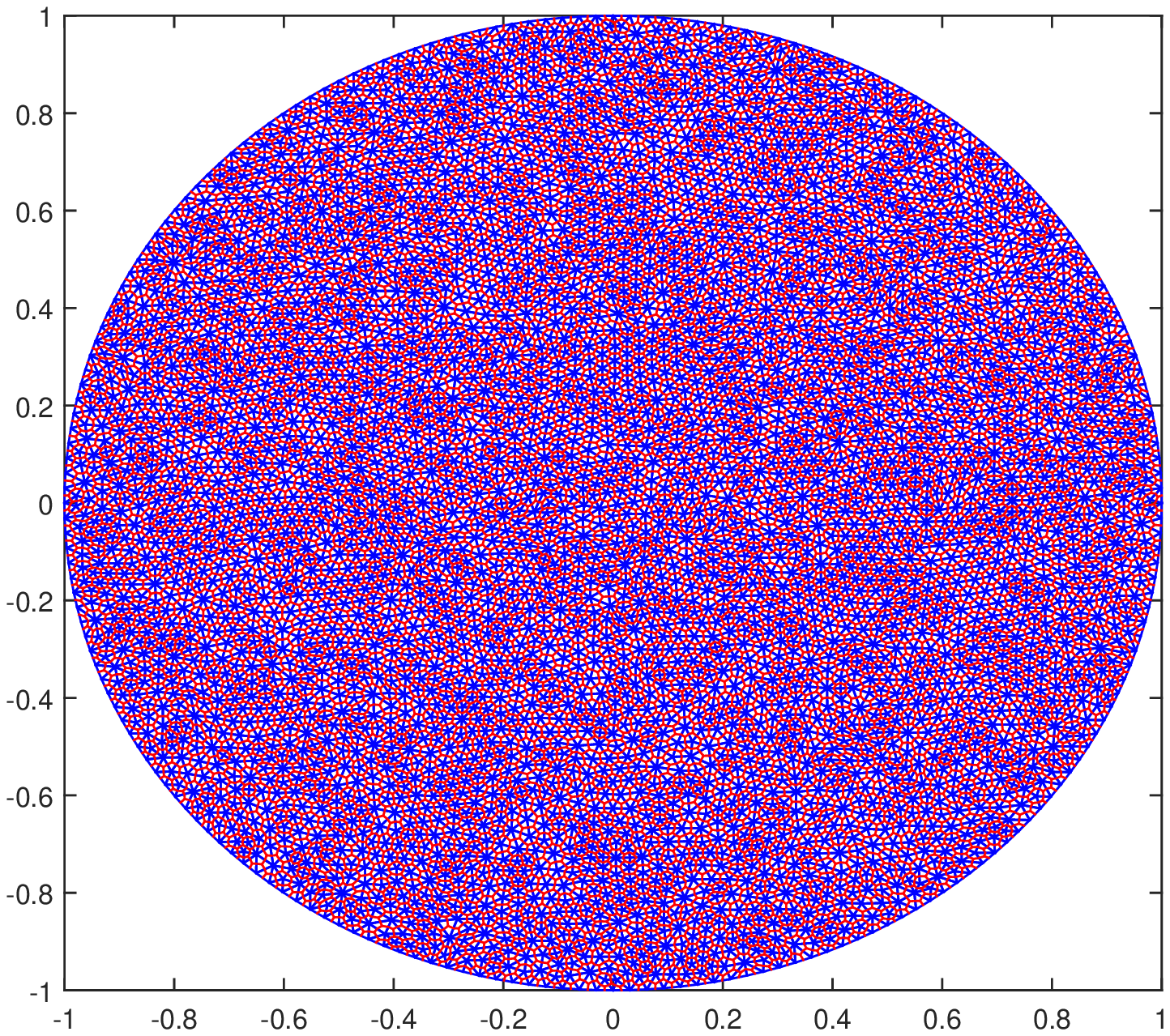}}
\caption{The unstructured meshes with control volumes for $h\approx  8.6550\times 10^{-2}, 4.5873\times 10^{-2}$, respectively}\label{fig9}
\end{center}
\end{figure}

The exact solution is given by $u(x,y,t)=e^{-t}({x^2}+{y^2}-1)^2$. Figure \ref{fig9} shows the circular domain partitioned by unstructured triangular meshes and control volumes for different $h$. In \cite{Yang17}, Yang et al. applied the Galerkin finite element method
for solving the two-dimensional Riesz space fractional diffusion equation with a nonlinear source term on convex domains. They developed an algorithm to form the stiffness matrix on triangular meshes, which can deal with space fractional derivatives on any convex domain. Here, we will make a comparison between our method (CVM) and Yang's method (FEM) for solving the two-dimensional Riesz space fractional diffusion equation (\ref{e21}) on a circular domain using the same triangular meshes. Firstly, we present a comparison of the density of the two stiffness matrices generated by FEM and CVM for different $h$ in Table \ref{tab6}. We can see that with $h$ decreasing the density of the two stiffness matrices reduces significantly. Compared to the stiffness matrix generated by FEM, the stiffness matrix generated by CVM is slightly more sparse. Next, we present a comparison of the error and convergence. Table \ref{tab7} displays the $L_{2}$ error, $L_{\infty}$ error and corresponding convergence order of $h$ for different $\alpha$, $\beta$ with $\tau=10^{-3}$ at $t=1$ by applying FEM.  Table \ref{tab8} highlights the error and convergence order by using FVM.  We can see that the accuracy of our method is similar to FEM, both of which are second order. Then, we present a comparison of CPU time for the two methods in Table \ref{tab9} both using the Bi-CGSTAB solver.  We choose $\alpha=\beta=0.8$ and $\tau=10^{-3}$ at $t=1$ to observe the running time for different $h$. We observe that compared to the running time of FEM, CVM can reduce the running time significantly, which illustrates that CVM is more effective for solving the two-dimensional Riesz space fractional diffusion equation on convex domains. This is mainly due to the bilinear form in \cite{Yang17} that involves 8 fractional derivative terms and the approximation of two-fold multiple integrals, which are approximated by Gauss quadrature, while for CVM we only need to calculate 4 fractional derivative terms and the approximation of line integrals. %In addition, we give a comparison of the exact solution $u(x,y,t)$ and numerical solution $u_h(x,y,t)$ in Figure \ref{fig10} and the error plot of $u(x,y,t)-u_h(x,y,t)$ in Figure \ref{fig11} for $h=4.5873\times 10^{-2}$, $\alpha=\beta=0.8$ with $\tau=10^{-3}$ at $t=1$ by applying CVM.
We can see that the numerical solution is in excellent agreement with the exact solution, which demonstrates the effectiveness of our numerical method again.

\begin{table}[h]
\begin{center}
\caption{The comparison of the density of stiffness matrix generated by FEM and CVM for different $h$}\label{tab6}
\begin{tabular}{ccccc}
\hline
\hline
  $N_e$     &     $h$    & Size           &  FEM          & CVM   \\
\hline
    174     & 2.8917E-01 &   $74\times74$   &  65.413 \%    &   55.332\%   \\
    570     & 1.6444E-01 &  $260\times260$  &  41.814 \%    &   33.521\%   \\
    2310    & 8.6550E-02 & $1104\times1104$ &  22.233 \%    &   17.469\% \\
    8744    & 4.5873E-02 & $4271\times4271$ &  11.712 \%    &   9.107\%\\
\hline
\hline
\end{tabular}
\end{center}
\end{table}

\begin{table}[h]
\begin{center}
\caption{The $L_{2}$ error, $L_{\infty}$ error and convergence order of $h$ for FEM with $\tau=10^{-3}$ at $t=1$}\label{tab7}
\begin{tabular}{cccccc}
\hline
\hline
 FEM &     $h$ & $L_{2}$ error & Order & $L_{\infty}$ error & Order \\
\hline
              & 2.8917E-01 & 6.7022E-03   &  --   & 5.8841E-03   & --     \\

$\alpha=0.80$ & 1.6444E-01 & 2.0787E-03   &  2.07 & 2.8557E-03   & 1.28    \\

$\beta=0.80$  & 8.6550E-02 & 5.2077E-04   &  2.16 & 8.1791E-04   & 1.95 \\

              & 4.5873E-02 & 1.3554E-04   &  2.12 & 2.3520E-04   & 1.96  \\
\hline
              & 2.8917E-01 & 6.9018E-03   & --    & 5.5925E-03  & --    \\

$\alpha=0.70$ & 1.6444E-01 & 2.1713E-03   & 2.05  & 2.7718E-03  & 1.24   \\

$\beta=0.90$  & 8.6550E-02 & 5.4452E-04   & 2.16  & 7.9048E-04  & 1.95  \\

              & 4.5873E-02 & 1.4147E-04   & 2.12  & 2.2242E-04  & 2.00  \\
\hline
\hline
\end{tabular}
\end{center}
\end{table}

\begin{table}[h]
\begin{center}
\caption{The $L_{2}$ error, $L_{\infty}$ error and convergence order of $h$ for CVM with $\tau=10^{-3}$ at $t=1$}\label{tab8}
\begin{tabular}{cccccc}
\hline
\hline
 CVM &     $h$ & $L_{2}$ error & Order & $L_{\infty}$ error & Order \\
\hline
              & 2.8917E-01 & 1.4782E-02   &  --   & 2.1786E-02  & --     \\

$\alpha=0.80$ & 1.6444E-01 & 4.5014E-03   & 2.11  & 7.5230E-03  & 1.88    \\

$\beta=0.80$  & 8.6550E-02 & 1.2275E-03   & 2.02  & 1.8279E-03  & 2.20  \\

              & 4.5873E-02 & 3.4069E-04   & 2.02  & 5.4557E-04  & 1.90  \\
\hline
              & 2.8917E-01 & 1.4950E-02   & --    & 2.1864E-02  & --     \\

$\alpha=0.70$ & 1.6444E-01 & 4.5530E-03   & 2.11  & 7.6462E-03  & 1.86    \\

$\beta=0.90$  & 8.6550E-02 & 1.2566E-03   & 2.01  & 1.8659E-03  & 2.20  \\

              & 4.5873E-02 & 3.4898E-04   & 2.02  & 5.4606E-04  & 1.94  \\
\hline
\hline
\end{tabular}
\end{center}
\end{table}

\begin{table}[H]
\begin{center}
\caption{The comparison of running time between FEM and CVM for different $h$ with $\alpha=\beta=0.80$, $\tau=10^{-3}$ at $t=1$}\label{tab9}
\begin{tabular}{cccc}
\hline
\hline
    $N_e$   &     $h$    &  FEM    & CVM   \\
\hline
    174     & 2.8917E-01 &  3.49 min    &   35.01 s   \\
    570     & 1.6444E-01 &  12.90 min   &   2.63min   \\
    2310    & 8.6550E-02 &  1.38 h      &   28.41min \\
    8744    & 4.5873E-02 &  17.89h      &   6.59h\\
\hline
\hline
\end{tabular}
\end{center}
\end{table}

\section{Conclusions}
In this paper, we considered the unstructured mesh control volume method for the two-dimensional space fractional diffusion equation with variable coefficients on convex domains. We partitioned the irregular convex domain using triangular meshes. Then we constructed the control volumes and solved the space fractional diffusion equation by utilising the finite volume method. Finally, numerical examples on irregular convex domains were studied, which verified the effectiveness and reliability of the method. We concluded that the numerical method can be extended to other arbitrarily shaped convex domains. Furthermore, according to the property of the stiffness matrix generated by the finite volume method, we chose a suitable sparse matrix format for the stiffness and utilised the Bi-CGSTAB iterative method to solve the linear system, which is more efficient than using Gauss elimination method. In addition, we made a comparison of our method with the finite element method proposed in \cite{Yang17}, which demonstrated that our method can reduce CPU time significantly while retaining the same accuracy and approximation property as the finite element method. In future work, we shall investigate the unstructured mesh control volume method applied to other fractional problems on irregular convex domains, such as the two-dimensional multi-term time-space fractional diffusion equation with variable coefficients or three-dimensional space fractional diffusion equations with variable coefficients.

%\section*{Acknowledgment}
%Authors Liu and Turner wish to acknowledge that this research was partially supported by the Australian Research Council (ARC) via the Discovery Project (DP180103858). Author Liu wishes to acknowledge that this research was partially supported by Natural Science Foundation of China (Grant No.11772046). Author Turner wishes to acknowledge that this research was also partially supported by the Australian Research Council (ARC) via the Discovery Project DP150103675.

\end{document}